\documentclass{amsart}
\usepackage[dvips]{graphicx}
\usepackage{amssymb,amsmath}
\newtheorem{theorem}{Theorem}[section]
\newtheorem{proposition}[theorem]{Proposition}
\newtheorem{lemma}[theorem]{Lemma}
\newtheorem{corollary}[theorem]{Corollary}

\theoremstyle{definition}

\theoremstyle{remark}
\newtheorem{remark}[theorem]{Remark}

\numberwithin{equation}{section}

\newcommand{\PSLC}{\mathrm{PSL}(2,\mathbb{C})}

\begin{document}

\title{A method to find ideal points from ideal triangulations}
\author{Yuichi KABAYA}
\address{Department of Mathematics, Tokyo Institute of Technology, 
2-12-1 Oh-okayama, Meguro-ku, Tokyo 152-8551, JAPAN}
\email{kabaya@math.titech.ac.jp}

\begin{abstract}
In Yoshida \cite{yoshida}, a method to find ideal points and boundary slopes from an ideal triangulation
was introduced.
But the method only gives a necessary condition which the valuation corresponding to an ideal point satisfies.
We give a simple sufficient condition for the existence of the corresponding ideal points.
\end{abstract}

\maketitle

\section{Introduction}
The character variety is an important object in 3-dimensional topology. 
In the seminal paper \cite{culler-shalen}, 
Culler and Shalen developed a relationship between incompressible surfaces and ideal points of a character variety.
Although Culler-Shalen theory gives much information about 3-manifolds,
it is difficult to compute the character variety.
On the other hand, it is easier to study only ideal points of a character variety.
In \cite{yoshida}, Yoshida introduced a method to find ideal points from an ideal triangulation of a 3-manifold.
His method gives (necessary) conditions which a valuation satisfies.
Since the conditions are given by linear equations, 
we can get candidates of ideal points only by solving linear equations.
But, by using this method, we can only find candidates of ideal points.
So we have to check that each of the candidates is actually corresponds to some ideal points.

Let $N$ be a compact orientable 3-manifold whose boundary is a torus.
Let $K$ be an ideal triangulation of $N$  with $n$ ideal tetrahedra and give a complex parameter 
$z_{\nu}$ for each ideal tetrahedron.
At each 1-simplex of $K$, we have the \emph{gluing equation} of $K$ as explained in
\cite{neumann-zagier}, \cite{yoshida}.
The gluing equations have the following form:
$
\displaystyle\prod_{\nu=1}^n z_{\nu}^{r'_{i,\nu}} (1-z_{\nu})^{r''_{i,\nu}} = \pm 1
$.
Let $\mathcal{D}(N)$ be the algebraic set defined by the gluing equations and 
call it the \emph{deformation variety}.
For each point of $\mathcal{D}(N)$ we can construct $\PSLC$-representation up to 
conjugacy.
So ideal points of $\mathcal{D}(N)$ are closely related to ideal points of the character variety
of $N$.
In order to find the ideal points of $\mathcal{D}(N)$, 
Yoshida \cite{yoshida} introduced a linear system of equations 
defined in terms of the above $(r'_{i,\nu},r''_{i,\nu})$.
Yoshida showed that in the case of the figure eight knot complement each integral solutions 
of the linear system of equations gives an ideal point. 
But in general this is not clear whether all integral solutions of the equations correspond to ideal points or not. 

In the paper under discussion, the author continues the work of Yoshida. 
He introduces a sufficient condition which guarantees that certain solutions of 
the linear system of equations corresponds actually to ideal points of $\mathcal{D}(N)$
(Theorem \ref{main-theorem}).
Let $I \in \{ 1,0,\infty\}^n$. $I$ describes how ideal tetrahedra of $K$ degenerate ($z_{\nu} \to 1$, $0$ or  $\infty$).
Then we introduce a vector $d(I)$ called \emph{degeneration vector} in Section \ref{sec:candidates}.
$d(I)$ can be calculated only by calculating determinants of some matrices with entry 
described in terms of $(r'_{i,\nu},r''_{i,\nu})$.
The key idea of the proof is that we can construct an appropriate embedding of $\mathcal{D}(N)$ 
into weighted projective space whose weight is given by the coefficients of degeneration vector.
This embedding gives a good compactification of $\mathcal{D}(N)$ near the ideal point corresponding to
the degeneration described by $I$.

This paper is organized as follows.
In section \ref{sec:basic} we explain some basic notion about ideal triangulations, character varieties,
and ideal points.
In section \ref{sec:candidates}, we introduce a method to find candidates of ideal points as explained in \cite{yoshida}.
In section \ref{sec:ideal}, we give a criterion for a candidate to be actually an ideal point.
In section \ref{sec:examples}, we give some examples from census manifolds.

\section{Basic definitions}
\label{sec:basic}
In this section, we review basic notions, ideal triangulations, character varieties, and ideal points.
References of ideal triangulation are  \cite{neumann-zagier} and \cite{neumann}.
For character varieties and ideal points, see \cite{culler-shalen}.
\subsection{Ideal triangulation}
Let $N$ be a compact oriented 3-manifold with the boundary $\partial N$ homeomorphic to a torus. 
Let $K$ be a cell complex which is given by gluing tetrahedra along faces.
We assume that the 0-simplex $K^{(0)}$ of $K$ is only one point $e^0$.
Let $Nbd(e^0)$ be a small neighborhood of $e^0$.
$K$ is an \emph{ideal triangulation} of $N$ 
if $K-Nbd(e^0)$ is homeomorphic to $N$.
We give an orientation induced by the orientation of $N$ for each tetrahedron.
Let $n$ be the number of ideal tetrahedra of $K$.
The number of 1-simplices of $K$ is $n$ since the Euler characteristic of $N$ vanishes.

\subsection{Ideal tetrahedron}
Let $\mathbb{H}^3$ be the upper half space model of the 3-dimensional hyperbolic space.
The ideal boundary of $\mathbb{H}^3$ is identified with $\mathbb{C}P^1$.
$\PSLC$ acts on $\mathbb{C}P^1$ by fractional linear transformations and the action extends to 
the isometry of $\mathbb{H}^3$.
An ideal tetrahedron is a geodesic 3-simplex with all vertices at $\mathbb{C}P^1$.
An ideal tetrahedron 
is described by 4 distinct points 
$z_0,z_1,z_2,z_3$ of $\mathbb{C}P^1$.
An edge of ideal tetrahedron is described by giving a pair of points $(z_i,z_j)$.
For each edge $(z_i,z_j)$, we define the complex parameter
$z=[z_i:z_j:z_k:z_l]=\frac{(z_k-z_j)(z_l-z_i)}{(z_k-z_i)(z_l-z_j)}$ where $k$ and $l$ are chosen in such a way that
$(i,j,k,l)$ forms the orientation of the ideal tetrahedron.
The complex number is not equal to 0 or 1 because $z_0,z_1,z_2,z_3$ are distinct.
Let $z$ be the complex parameter of an edge.
Then the opposite edge has the same complex parameter.
For the other edges, the complex parameters are given by $z'_i=\frac{1}{1-z_i}$, $z''_i=1-\frac{1}{z_i}$.
We put $w_i=1-z_i$ then we have $z'_i=1/w_i$ and $z''_i=-w_i/z_i$.

Let $K$ be an ideal triangulation of $N$.
Let $e_i$ be a 1-simplex of $K$.
There are the edges of ideal tetrahedra attached to $e_i$ each of which has complex parameter
$z_\nu$, $\frac{1}{1-z_\nu}$ or $1-\frac{1}{z_\nu}$.
Let $p_{i,\nu}$ be the number of edges attached to $e_i$ whose complex parameters are $z_{\nu}$.
(So $0 \leq p_{i,\nu} \leq 2$.)
We also define $p'_{i,\nu}$ and $p''_{i,\nu}$ be the number of edges attached to $e_i$ whose complex parameters 
are $\frac{1}{1-z_{\nu}}$ and $1-\frac{1}{z_\nu}$ respectively.
Then we define $R_i$ by
\[
R_i=\prod_{\nu=1}^n (z_{\nu})^{p_{i,\nu}} (z'_{\nu})^{p_{i,\nu}'} (z''_{\nu})^{p_{i,\nu}''} \quad (i=1, \dots ,n).
\]
$R_i$ is simplified as follows:
\[
R_i  = \prod_{\nu=1}^n (-1)^{p''_{i,\nu}} z_{\nu}^{(p_{i,\nu}-p''_{i,\nu})}
w_{\nu}^{(p''_{i,\nu}-p'_{i,\nu})} 
 =\prod_{\nu=1}^n (-1)^{p''_{i,\nu}} z_{\nu}^{r'_{i,\nu}} w_{\nu}^{r''_{i,\nu}}
\]
where
\[
r'_{i,\nu}=p_{i,\nu}-p''_{i,\nu}, \quad r''_{i,\nu}=p''_{i,\nu}-p'_{i,\nu}.
\]
We remark that each $(z_1, \dots, z_n)$ satisfying $R_i=1$  $(i=1, \dots, n)$ gives a representation of $\pi_1(N)$ into 
$\PSLC$ (see subsection \ref{deformation}).
Because $R_1 \dots R_n = 1$, we only have to consider $n-1$ equations.
So we omit n-th equation from the edge relations.
We define $r_i=(r'_{i,1},r''_{i,1}, \dots, r'_{i,n},r''_{i,n})$ and 
\begin{equation}
\label{R}
R=\begin{pmatrix} r_1 \\ \vdots \\ r_{n-1} \end{pmatrix}
=\begin{pmatrix} 
r_{1,1} & \hdots & r_{1,n} \\
\vdots &         & \vdots \\
r_{n-1,1} & \hdots & r_{n-1,n}
\end{pmatrix}.
\end{equation}


Let $\mathcal{M}$ and $\mathcal{L}$ be simple closed curves on $\partial N$ 
which generate $H_1(\partial N,\mathbb{Z})$.
We define a pair of integers $(m_{\nu}',m_{\nu}'')$ for $\mathcal{M}$ in the following manner
(also define $(l_{\nu}',l_{\nu}'')$ for $\mathcal{L}$).
When $\mathcal{M}$ passes through the boundary torus $\partial N$, we homotope $\mathcal{M}$ 
so that $\mathcal{M}$ does not meet
0-simplices of $\partial N$ and passes through each 2-simplex of $\partial N$ from one edge to another edge.
As $\mathcal{M}$ passes through a 2-simplex of $\partial N$, 
we assign a complex parameter corresponding to the vertex binding these two edges.
Let $M$ be the multiplication of these complex parameters or inverses of them  
according as  $\mathcal{M}$ passes the tetrahedron anti-clockwise or clockwise viewing from 
the vertex of $K$.
Then $M$ (and $L$) can be written as 
\[
M=\pm\prod_{\nu=1}^n z_{\nu}^{m_{\nu}'}(1-z_{\nu})^{m_{\nu}''}, \quad 
L=\pm\prod_{\nu=1}^n z_{\nu}^{l_{\nu}'}(1-z_{\nu})^{l_{\nu}''} ,
\]
where $(m'_{\nu},m''_{\nu})$ (and  $(l'_{\nu},l''_{\nu})$) are some integers.
When $(z_1,\dots, z_n)$ satisfies $R_i=1$ $(i=1,\dots,n-1)$,
these represent the derivatives of holonomies along $\mathcal{M}$ and $\mathcal{L}$
(see Neumann-Zagier \cite{neumann-zagier}).
We define $m=(m'_1,m''_1 \dots, m'_n,m''_n)$ and $l=(l'_1,l''_1, \dots, l'_n,l''_n)$.

$\mathbb{R}^{2n}$ has a symplectic form given by 
\[
x \wedge y = \sum_{k=1}^{n} x'_k y''_k - x''_k y'_k
\]
for $x=(x'_1, \dots, x'_n, x''_1, \dots x''_n)$ and 
$y=(y'_1, \dots, y'_n, y''_1, \dots y''_n)$. 
Let $[R]= \mathrm{span}_{\mathbb{R}} \langle r_1, \dots, r_{n-1}\rangle$.
We denote the orthogonal complement of $[R]$ with respect to $\wedge$ by $[R]^{\bot}$. 

\subsection{$\PSLC$ representations}
Let $R(N)=\mathrm{Hom}(\pi_1(N),\PSLC)$ be the $\PSLC$-representation variety of $N$.
In this paper, we use the term \emph{variety} for an algebraic set so a variety may not be irreducible.
Since $\PSLC$ is an affine algebraic group, $R(N)$ is an affine algebraic set.
$\PSLC$ acts on $R(N)$ by conjugation and also acts on the coordinate ring of $R(N)$.
We denote the coordinate ring of $R(N)$ by $\mathbb{C}[R(N)]$.
The character variety $X(N)$ of $N$ is the affine variety of $\PSLC$-invariant
subring $\mathbb{C}[R(N)]^{\PSLC}$.
It is known that there is a surjective regular map $t:R(N) \to X(N)$. 
$X(N)$ can be regarded as the set of the squares of characters \cite{heusener-porti}. 
For any algebraic curve $Y \subset X(N)$, Boyer and Zhang \cite{boyer-zhang} constructed a Culler-Shalen theory:
 there is an incompressible surface for each ideal point of $Y$.

\subsection{Deformation variety}
\label{deformation}
Let
\begin{equation}
\label{deformation-curve}
\begin{split}
\mathcal{D}(N,K)  & =\{(z_1, \dots, z_n, w_1, \dots, w_n) \in (\mathbb{C}^*)^{2n} | \\
& R_1(z,w)=1, \dots, R_{n-1}(z,w)=1, z_1+w_1=1, \dots, z_n+w_n=1 \} \\
    & = \{(z_1, \dots, z_n) \in (\mathbb{C}-\{0,1\})^n | R_1(z)=1, \dots, R_{n-1}(z)=1 \}.
\end{split}
\end{equation}
We call $\mathcal{D}(N,K)$ the \emph{deformation variety} and the equations $R_i=1$ \emph{gluing equations}.
In this paper, we assume that $\mathcal{D}(N,K)$ is not empty.
We often simply denote $\mathcal{D}(N,K)$ by $\mathcal{D}(N)$. 
$\mathcal{D}(N)$ can be represented as an affine algebraic variety as follows:
\[
\begin{split}
\mathcal{D}(N,K) & =
\{ (z_1, \dots, z_n, w_1, \dots, w_n,t) \in \mathbb{C}^{2n+1} | R_1(z,w)=1, \dots, R_{n-1}(z,w)=1, \\
& z_1+w_1=1, \dots, z_n+w_n=1 ,  z_1 \cdots z_n w_1 \cdots w_n t =1 \} \subset \mathbb{C}^{2n+1}.
\end{split}
\]
In (\ref{deformation-curve}), $\mathcal{D}(N,K)$ is represented as the solutions of $n-1$ equations in $(\mathbb{C}-\{0,1\})^n$,  
the dimension of $\mathcal{D}(N)$ is equal to or greater than 1.
So $\mathcal{D}(N)$ contains an algebraic curve.

The following proposition is well-known (see \cite[p163-164]{yoshida}).
\begin{proposition}
For each point $p \in \mathcal{D}(N)$,  we can construct $\PSLC$ representation of $\pi_1(N)$.
This defines an algebraic map $\mathcal{D}(N) \to X(N)$.
\end{proposition}
Roughly speaking,
a point $p \in \mathcal{D}(N)$ defines the developing map, and the developing map gives the holonomy map $\pi_1(N) \to \PSLC$ which
is invariant under conjugation.
Since this construction is algebraic, so we obtain the algebraic map $\mathcal{D}(N) \to X(N)$.

\subsection{Ideal points of an affine variety and valuations}
\label{ideal-points-of-affine}
In this subsection we define the ideal points of an affine algebraic curve and give an exposition of valuation theory.
For details, see the original paper of Culler and Shalen \cite{culler-shalen}.

Let $C$ be an affine algebraic curve.
Let $\bar{C}$ be a projective completion and let $\tilde{C}$ be the projective smooth curve which is birational equivalent to $C$.
Then there is a birational map $f : \tilde{C} \to \bar{C}$. 
Since $\tilde{C}$ is smooth, $f$ is regular.
The points $f^{-1}(\bar{C}-C) \subset \tilde{C}$ are called the \emph{ideal points} of $C$.

Let $\tilde{C}$ be a smooth projective curve and $\mathbb{C}(\tilde{C})$ be the function field of $\tilde{C}$.
A \emph{valuations} of $\mathbb{C}(\tilde{C})/\mathbb{C}$ is a map $v: \mathbb{C}(\tilde{C})-\{0\} \to \mathbb{Z}$ satisfying
(1) $v(xy)=v(x)+v(y)$, (2) $v(x+y) \geq \min(v(x),v(y))$ for all $x,y \in \mathbb{C}(\tilde{C})$ and 
(3) $v(x)=0$ for $x \in \mathbb{C}$. 
We simply call valuations of $\mathbb{C}(\tilde{C})/\mathbb{C}$ by valuation of $\mathbb{C}(\tilde{C})$.
There is a one-to-one correspondence between the valuations of $\mathbb{C}(\tilde{C})/\mathbb{C}$ and the points of $\tilde{C}$. 
Let $C$ be an affine algebraic curve which is birational equivalent to $\tilde{C}$.
Let $p$ be a point of $C$, then the corresponding valuation $v$ satisfies $v(f) \geq 0$ for any regular function $f \in \mathbb{C}[C]$.

We return to our problem.
Since $\mathcal{D}(N)$ is an affine algebraic variety, the notion of ideal points of $\mathcal{D}(N)$ is well-defined.
Roughly speaking, the ideal point of $\mathcal{D}(N)$ is the point whose parameters
$z_i$ are equal to $0$, $\infty$ or $1$. 
Fix an irreducible curve $Y_0$ of $\mathcal{D}(N)$.
Let $X_0$ be the closure of the image of $\mathcal{D}(N) \to X(N)$.
Then we have an extension of field $\mathbb{C}(X_0) \subset \mathbb{C}(Y_0)$.
Let $p$ be an ideal point of $Y_0$ and $v$ be the corresponding valuation.
$v$ also defines a valuation on $\mathbb{C}(X_0)$ by restriction.
We remark that $v$ may not correspond to an ideal point of $\mathbb{C}(X(N))$
i.e. $v$ may satisfy $v(f) \geq 0$ for any regular function $f \in \mathbb{C}[X_0]$.

\section{Candidates of ideal points}
\label{sec:candidates}
Let $p$ be an ideal point of $\mathcal{D}(N)$ and $v$ be the associated valuation.
Then $v$ satisfies
\[
\begin{split}
0 & = v(1)=v(R_i) = \displaystyle\sum_{\nu=1}^{n} \biggl( r'_{i,\nu} v(z_\nu ) + r''_{i,\nu} v( w_\nu ) \biggr) \\
  &  = (r'_1,r''_1, \dots, r'_n,r''_n) \wedge (-v(w_1), v(z_1), \dots, -v(w_n), v(z_n)).
\end{split}
\]
For the valuation $v$, there is a sequence of points $(z(k)_1,w(k)_1,\dots,z(k)_n,w(k)_n) \in \mathcal{D}(N)$, $k=1, 2,\dots$, so that
\[
\frac{(\log |z(k)_1|,\log |w(k)_1|, \dots,\log |z(k)_n|,\log |w(k)_n|)}
{(\sqrt{1+\sum_{\nu=1}^{n}(\log |z(k)_\nu|^2+\log |w(k)_\nu|^2)})}
\]
converges to $(-v(z_1),-v(w_1), \dots , -v(z_n),-v(w_n))$
(see Tillmann \cite{tillmann-ideal}).
Because $z_\nu+w_\nu=1$, $(-\log |w_\nu |,\log|z_\nu|)$ diverge to the direction
 $(1,0)$, $(0,-1)$ or $(-1,1)$ when $z_\nu$ converges to $1$, $0$ or $\infty$ respectively.
Let $\rho_1=(1,0)$, $\rho_0=(0,-1)$ and $\rho_{\infty}=(-1,1)$.
The valuation  $v$ satisfies 
\[
(-v(w_1), v(z_1), \dots -v(w_n),v(w_n))= (t_1 \rho_{i_1}, \dots, t_n \rho_{i_n}) \quad (t_i \geq 0)
\]
and 
\[
r_i \wedge (-v(w_1), v(z_1), \dots -v(w_n),v(w_n))=0 \quad (i=1, \dots, n-1).
\]
In this section, we show how to compute the solutions of these linear equations.

\subsection{}
Let $I=(i_1, \dots , i_n) \in \{ 1,0, \infty\}^n$ and call it a \emph{degeneration index}.
A degeneration index describes degeneration of ideal tetrahedra.
Let
\[
\begin{split}
S(I)=S(i_1,i_2, \dots, i_n)=\{ (t_1 \rho_{i_1}, \dots, t_n \rho_{i_n}) | t_i \in \mathbb{R} \} 
\subset \mathbb{R}^{2n}, \\
H(I)=H(i_1,i_2, \dots, i_n)=\{ (t_1 \rho_{i_1}, \dots, t_n \rho_{i_n}) | t_i \in \mathbb{R}, t_i \geq  0 \} 
\subset \mathbb{R}^{2n}
\end{split}
\]
for a degeneration index $I$.
$S(I)$ is an n-dimensional subspace of $\mathbb{R}^{2n}$ and $H(I)$
is a convex cone in $S(I)$.
In this subsection we will show how to describe and compute the set
\begin{equation}
\label{necessary condition}
\bigcup_{i_1, \dots, i_n} H(i_1, \dots, i_n) \bigcap [R]^{\bot}.
\end{equation}
As explained in the beginning of this section, this set can be regarded as candidates of ideal points 
\cite{yoshida}. 
This set also describes normal surfaces in $N$ \cite{tillmann-normal}.

Fix a degeneration index $I=(i_1, \dots, i_n)$.
Because $H(I) \bigcap [R]^{\bot}$ is a convex cone in $\mathbb{R}^{2n}$, 
we want to find a set of generators of this convex cone.
In general, $S(I) \bigcap [R]^{\bot}$ is a 1-dimensional subspace because
$S(I)$ and $[R]^{\bot}$ have dimension $n$ and $n-1$ respectively.
So $H(I) \bigcap [R]^{\bot}$ is a 1-dimensional convex cone or the origin in general.
To determine the exact dimension of $H(I) \bigcap [R]^{\bot}$, we introduce some definitions.
Let
\[
r(I)_{j,\nu}=
\left\{ \begin{array}{ll}
r''_{j,\nu}  & \textrm{if $i_\nu = 1$}\\
r'_{j,\nu}  & \textrm{if $i_\nu = 0$ }\\
-r'_{j,\nu}-r''_{j,\nu}  & \textrm{if $i_\nu = \infty$}.
\end{array} \right. 
\]

We define a \emph{degeneration matrix} of $I$ by
\[
R(I)=
\begin{pmatrix} 
r(I)_{1,1} & \hdots & r(I)_{1,n} \\
\vdots & & \vdots \\
r(I)_{n-1,1} & \hdots & r(I)_{n-1,n}
\end{pmatrix}.
\]
\begin{lemma}
\label{dimension}
$\ker (R(I))$ and $S(I) \cap [R]^{\bot}$ are isomorphic by the morphism
$\ker (R(I)) \to S(I) \cap [R]^{\bot} : (a_1, \dots, a_n) \mapsto (a_1 \rho_{i_1}, \dots, a_n \rho_{i_n})$.
Therefore the dimension of the subspace $ S(I) \bigcap [R]^{\bot}$ is equal to $n-\mathrm{rank}(R(I))$.
\end{lemma}
\begin{proof}
We have 
\[
a_\nu \rho_{i_\nu} \wedge (r'_{j,\nu},r''_{j,\nu})= 
\left\{ \begin{array}{ll}
a_\nu r''_{j,\nu}  & \textrm{if $i_\nu = 1$}\\
a_\nu r'_{j,\nu}  & \textrm{if $i_\nu = 0$ }\\
a_\nu(-r'_{j,\nu}-r''_{j,\nu})  & \textrm{if $i_\nu = \infty$}.
\end{array} \right. 
\]
where $r_i=(r'_{i,1},r''_{i,1}, \dots, r'_{i,n},r''_{i,n})$.
Therefore 
\[
R(I) \begin{pmatrix} a_1 \\ \vdots \\ a_n \end{pmatrix} 
= \begin{pmatrix} 
(a_1 \rho_{i_1}, \dots, a_n \rho_{i_n}) \wedge r_1 \\ \vdots \\ 
(a_1 \rho_{i_1}, \dots, a_n \rho_{i_n}) \wedge r_{n-1} 
\end{pmatrix}.
\]
This means that the map $\ker (R[I]) \to S(I) \cap [R]^{\bot}$ is well-defined.
We can also construct the inverse map, so the map is an isomorphism.
\end{proof}
We define
\[
d(I)_{\nu}= (-1)^{\nu+1} \mathrm{det} 
\begin{pmatrix} 
r(I)_{1,1} & \hdots & \widehat{r(I)_{1,\nu}}   & \hdots & r(I)_{1,n} \\
\vdots   &        & \vdots     &        & \vdots \\
r(I)_{n-1,1} & \hdots & \widehat{r(I)_{n-1,\nu}}   & \hdots       & r(I)_{n-1,n}
\end{pmatrix}
\]
where the hat means removing the column.
Then we define a \emph{degeneration vector} by
\[
d(I)=(d(I)_1,d(I)_2, \dots, d(I)_n) \in \mathbb{Z}^n \subset \mathbb{R}^n.
\]
\begin{lemma}
\label{determinant}
A degeneration vector corresponds to an element of $ S(I) \cap [R]^{\bot}$.
If $ S(I) \cap [R]^{\bot}$ is one dimensional, then this vector is a generator of 
$ S(I) \cap [R]^{\bot}$
(otherwise this vector is zero).
\end{lemma}
\begin{proof}
We have
\[
\begin{vmatrix} 
r(I)_{k,1} & \dots & r(I)_{k,n} \\ 
r(I)_{1,1} & \dots & r(I)_{1,1} \\
\vdots    &       & \vdots \\
r(I)_{n-1,1} &   & r(I)_{n-1,n}
\end{vmatrix} = 0
\]
for $k = 1, \dots ,n-1$.
So we obtain $r(I)_{k,1}d(I)_1 + \dots + r(I)_{k,n} d(I)_n$ for all $k = 1, \dots, n$.
This means that $d(I) \in \ker(R(I))$.
\end{proof}

We denote $d(I)>0$ ($d(I) \geq 0$) if $d(I)_{\nu} >  0$ ($d(I)_{\nu} \geq 0$) for all $\nu$.
We denote $d(I)<0$ and $d(I) \leq 0$ similarly.  
If $d(I) \geq 0$ or $d(I) \leq 0$, we obtain an element of $ H(I) \bigcap [R]^{\bot}$
by multiplication by $-1$ if necessary.

\subsection{}
\label{compute-cone}
In case that the dimension of $S(I) \cap [R]^{\bot}$ is greater than 1, we also have an algorithm 
to compute generators of the cone $H(I) \cap [R]^{\bot}$.
We define
\[
\begin{split}
S(I)(\epsilon_1, \dots, \epsilon_n )
=\{ ((\epsilon_1 a_1) \rho_{i_1}, \dots, (\epsilon_n a_n) \rho_{i_n}) | a_i \in \mathbb{R} \}, \\
H(I)(\epsilon_1, \dots, \epsilon_n )
=\{ ((\epsilon_1 a_1) \rho_{i_1}, \dots, (\epsilon_n a_n) \rho_{i_n}) | a_i \in \mathbb{R}, a_i \geq  0 \} 
\end{split}
\]
for $\epsilon = (\epsilon_1, \dots, \epsilon_n) \in \{ 0,1\}^n$.
$H(I)(\epsilon)$ is a face of the cone $H(I)$.
We define a matrix $R(I)(\epsilon)$ by omitting columns corresponding to $\epsilon_\nu=0$.
As in Lemma~\ref{dimension}, $S(I)(\epsilon) \cap [R]^{\bot}$ is isomorphic to the kernel of $R(I)(\epsilon)$.
Then the generators of $S(I) \cap [R]^{\bot}$ can be computed similar to Lemma \ref{determinant} as follows.
By multiplying elementary matrices from left, permuting the indices of ideal simplices and removing trivial rows,
we can assume that $R(I)$ is equal to 
$\begin{pmatrix} 
1 &    &     &  c_{1,1} & \hdots & c_{1,n-r} \\
 &  \ddots &   &   \vdots &        & \vdots \\
 & & 1       &  c_{r,1} & \hdots & c_{r,n-r}  \\
\end{pmatrix}$
where $r$ $( < n-1)$ is the rank of $R(I)$.
Consider $\epsilon=(\epsilon_1, \dots, \epsilon_n)$ with $n-1-r$ entries satisfying $\epsilon_\nu=0$. 
Remove the columns with $\epsilon_\nu=0$ from $R(I)$.
Then we have an $r \times (r+1)$ matrix for each $\epsilon$.
Then we can also define a degeneration vector for this matrix.
If all the coefficients of the degeneration vector have same sign, this vector corresponds to a generator of 
$H(I) \cap [R]^{\bot}$.

\section{Ideal points}
\label{sec:ideal}
\subsection{Ideal points of $\mathcal{D}(N)$}
In the previous section, we observed that there is a necessary condition that
the valuation corresponding to an ideal point must satisfy (\ref{necessary condition}).
In this section we give a criterion which guarantees that 
the candidate of ideal point actually corresponds to an ideal point.
\begin{theorem}
\label{main-theorem}
Let $I=(i_1, \dots, i_n)$ be an element of $\{ 1,0,\infty \}^n$.
If $d(I)>0$ or $d(I)<0$
then there are ideal points of $\mathcal{D}(N)$ corresponding to $I$. 
The number of the ideal points is $\gcd(d(I)_1, \dots , d(I)_n)$.
\end{theorem}
We remark that the ideal points of $\mathcal{D}(N)$ may not correspond to the ideal points of $X(N)$.
We will study the ideal points of $X(N)$ at the next subsection \ref{psl-ideal}.

We will prove this theorem by embedding $\mathcal{D}(N)$ into a weighted projective space and studying the local behavior
of the points at infinity. 

A \emph{weighted projective space} with index $(m_0, \dots, m_n) \in \mathbb{Z}^{n+1}$ is
the quotient of $\mathbb{C}^{n+1} - \{(0,\dots,0) \} $ by the following equivalence relation
\[
(z_0, \dots ,z_n) \sim (c^{m_0} z_0 , \dots , c^{m_n} z_n) \textrm{ for some $c \in \mathbb{C}^*$} 
\]
and denoted by $\mathbb{C}P(m_0,\dots,m_n)$.
We denote an equivalent class by $[z_0; \dots ;z_n]$.
For instance, the weighted projective space with index $(1, \dots, 1)$ is the projective space.
Let $U_i = \{ [z_0; \dots; z_{i-1};1;z_{i+1}; \dots ;z_n]\} \subset \mathbb{C}P(m_0,\dots,m_n)$.
A weighted projective space $\mathbb{C}P(m_0, \dots ,m_n)$ is covered by $n+1$ charts $U_i$. 
Let $\mathbb{C}^{n} \to U_i$ be the map defined by
$(z_1, \dots ,z_n) \mapsto [z_1; \dots ;z_{i-1};1;z_{i+1}; \dots ;z_n]$.
The $m_i$-th root of unity $\zeta_{m_i}$ acts on $\mathbb{C}^n$ by 
$(z_1, \dots, z_n) \mapsto ( (\zeta_{m_i})^{m_0} z_1, \dots, 1, \dots, (\zeta_{m_i})^{m_n} z_n )$.
The quotient by this action induces a biholomorphic map $\mathbb{C}^{n}/\langle\zeta_{m_i}\rangle  \to U_i$.
We call the map $\mathbb{C}^{n}  \to U_i$ \emph{inhomogeneous coordinate}
although this is not a coordinate in a usual sense.
We have to factor out by the action of the cyclic group $\langle \zeta_{m_i} \rangle$.
Let $n(z_1,\dots,z_n)=\gcd( \{ m_l | z_l \neq 0 \} )$.
The following lemma is obvious by the definition.
\begin{lemma}
\label{freeness}
The action of $\zeta_{m_i}$ on $(z_1, \dots , z_n) \in \mathbb{C}^n$ is free if and only if 
$n(z_1,\dots,z_n)$ and $m_i$ are relatively prime.
\end{lemma}
In particular, if $m_i=\pm 1$,
$(\mathbb{C}^*)^n$ is biholomorphically embedded in $\mathbb{C}P(m_0,\dots,m_n)$.

\begin{proof}[proof of Theorem \ref{main-theorem}]
We fix $I$ which satisfies $d(I)>0$ or $d(I)<0$.
If $d(I)<0$, we multiply $-1$ so that $d(I)>0$.
Let $c=\gcd(d(I)_1,\dots,d(I)_n)$ and $d'(I)_\nu=d(I)_\nu /c$.
We often abbreviate $d(I)_\nu$ and $d'(I)_\nu$ to $d_\nu$ and $d'_\nu$.
We substitute 
\begin{equation}
\label{coordinate}
\begin{split}
w_{\nu}=a_{\nu} t^{d'_{\nu}} \quad \mathrm{if} \quad i_{\nu} = 1, \\
z_{\nu}=a_{\nu} t^{d'_{\nu}} \quad \mathrm{if} \quad i_{\nu} = 0, \\
1/{z_{\nu}} = a_{\nu} t^{d'_{\nu}} \quad \mathrm{if} \quad i_{\nu} = \infty .
\end{split}
\end{equation}
We remark that $z_{\nu} \to 1, 0, \infty$ for $i_\nu=1,0,\infty$ respectively as $t \to 0$. 
We define 
\[
\overline{r(i_\nu)}_{j,\nu}=\overline{r(I)}_{j,\nu}=
\left\{ 
\begin{array}{ll}
r'_{j,\nu}  & \textrm{if $i_\nu = 1$}\\
r''_{j,\nu}  & \textrm{if $i_\nu = 0$ }\\
r''_{j,\nu}  & \textrm{if $i_\nu = \infty$}.
\end{array} \right. 
\]
then
\[
\begin{split}
z_{\nu}^{r'_{i,\nu}} w_{\nu}^{r''_{i,\nu}} &=
\left\{ \begin{array}{ll}
(w_{\nu})^{r(I)_{i,\nu}} (z_{\nu})^{\overline{r(I)}_{i,\nu}}  &  \textrm{if $i_\nu = 1$} \\
(z_{\nu})^{r(I)_{i,\nu}} (w_{\nu})^{\overline{r(I)}_{i,\nu}}  &  \textrm{if $i_\nu = 0$} \\
(z_{\nu}^{-1})^{r(I)_{i,\nu}} (z_{\nu}^{-1} w_{\nu})^{\overline{r(I)}_{i,\nu}}  &  \textrm{if $i_\nu = \infty$} 
\end{array} \right. \\
& =\pm a_{\nu}^{r(I)_{i,\nu}} (1-a_{\nu} t^{d'_{\nu}})^{\overline{r(I)}_{i,\nu}} .
\end{split}
\]
Then the gluing equations $R_i=1$ are replaced by
\begin{equation}
\label{replaced}
R_i(t,a_1,\dots ,a_n)
=\pm\displaystyle\prod_{\nu =1}^{n} a_{\nu}^{r(I)_{i,\nu}} (1-a_{\nu} t^{d'_{\nu}})^{\overline{r(I)}_{i,\nu}}= 1
\quad (i=1, \dots n-1).
\end{equation}
The system of equations (\ref{replaced}) is well-defined on $\mathbb{C}P(-1,d'_1,\dots,d'_n)$.
We abbreviate $\mathbb{C}P(-1,d'_1,\dots,d'_n)$ to $\mathbb{C}P(-1,d')$.
The map $\phi: (\mathbb{C}^*)^n \to \mathbb{C}P(-1,d'_1, \dots,d'_n)$
defined by $(z_1, \dots , z_n) \mapsto [t;a_1; \dots; a_n]=[1; a_1 t^{d'_1}; \dots; a_n t^{d'_n}]$ 
is a biholomorphic embedding by Lemma \ref{freeness} 
where $a_i$ and $t$ are defined in (\ref{coordinate}).
Let $\widetilde{\mathcal{D}(N)}$ be the set satisfying the equations (\ref{replaced}).
Since $\phi(\mathcal{D}(N))$ satisfies the equations (\ref{replaced}), $\mathcal{D}(N)$ is mapped into $\widetilde{\mathcal{D}(N)}$.
We call $\mathbb{C}P(d')=\{ [0;a_1; \dots; a_n]\} \subset \mathbb{C}P(-1,d')$ 
the hyperplane at infinity.
The set $\widetilde{\mathcal{D}(N)}-\phi(\mathcal{D}(N)) \supset \widetilde{\mathcal{D}(N)} \cap \mathbb{C}P(d')$ 
contains ideal points of $\mathcal{D}(N)$.
We have
\[
\begin{split}
\widetilde{\mathcal{D}(N)} \cap \mathbb{C}P(d') = \{ (a_1, \dots, a_n) | R_i(0,a_1,\dots, a_n)=1 \quad (i=1,\dots n-1)\} \\
= \{ (a_1, \dots, a_n) |
\displaystyle\prod_{\nu=1}^{n} a_{\nu}^{r(I)_{i,\nu} } = \pm 1  \quad (i=1, \dots, n-1) \}.
\end{split}
\]
So we study the system of equations:
\begin{equation}
\label{eq:infinity}
\displaystyle\prod_{\nu=1}^{n} a_{\nu}^{r(I)_{i,\nu} } = \pm 1 \quad (i=1, \dots, n-1).
\end{equation}
Even if we replace one relation $R_i=1$ by $R_i R_j^{n}=1$ where $j \neq i$ and $n \in \mathbb{Z}$,
the set of solutions $\widetilde{\mathcal{D}(N)} \cap \mathbb{C}P(d')$ does not change.
If we describe the system of equations by the matrix
\[
\displaystyle\prod_{\nu=1}^{n} a_{\nu}^{r(I)_{i,\nu} } = \pm 1 \quad (i=1, \dots, n-1)
 \longleftrightarrow 
\begin{pmatrix}
r(I)_{1,1} & \hdots & r(I)_{1,n} \\
\vdots  &  &  \\
r(I)_{n-1,1} & \hdots  &  r(I)_{n-1,n} 
\end{pmatrix},
\]
the above operation corresponds to adding some integer multiple of the j-th row to the i-th row.
We can reduce the matrix to the following form
\[
\begin{pmatrix}
c_{1,1} &  c_{1,2}  &  \hdots   & c_{1,n-1}   &  e_1 \\
        &  c_{2,2}  &           & \vdots    &  \vdots  \\
        &           &  \ddots   &  \vdots   &  \vdots \\
        &           &           & c_{n-1,n-1} &  e_{n-1}
\end{pmatrix}.
\]
In these operations, the degeneration vector does not change.
In particular, $d_n=(-1)^n c_{1,1} c_{2,2} \dots c_{n-1,n-1}$.
By our assumption that $d_n \neq 0$, we have $c_{i,i} \neq 0$ for all $i$.
Since $d_{n-1} \neq 0$, we also have $e_{n-1} \neq 0$.
The last row of the above matrix corresponds to the equation $a_{n-1}^{c_{n-1,n-1}} a_{n}^{e_{n-1}}=\pm 1$
and therefore $a_{n-1} \neq 0$ and $a_n \neq 0$.
By the $(n-2)$-th row, we can show that $a_{n-2} \neq 0$ as before.
Inductively, we can conclude that $a_\nu \neq 0$ for all $\nu$.
The equations (\ref{eq:infinity}) have an ambiguity arising from the weighted multiplication of $\mathbb{C}^*$,
we solve the equations in the inhomogeneous coordinate $U_n$. 
The equations (\ref{eq:infinity}) in $U_n$ are obtained by substituting $a_n=1$.  
So the resulting equations are 
\begin{equation}
\label{at-infinity}
\begin{split}
a_1^{c_{1,1}} a_2^{c_{1,2}} \dots a_{n-1}^{c_{1,n}} & = \pm 1, \\
\quad \quad \quad \vdots \\
a_{n-1}^{c_{n-2,n-2}} a_{n-1}^{c_{n-2,n-1}} & = \pm 1, \\
a_{n-1}^{c_{n-1,n-1}} & = \pm 1.
\end{split}
\end{equation}
The solution in $U_n$ is obtained by factoring the solutions of (\ref{at-infinity}) 
by the action of $\zeta_{|d'_n|}$.
By the $(n-1)$-th equation, $a_{n-1}$ has $|c_{n-1,n-1}|$ solutions.
Substitute each solution $a_{n-1}$, 
$a_{n-2}$ has $|c_{n-2,n-2}|$ solutions by the $(n-2)$-th equation.
So $(a_{n-2},a_{n-1})$ has $|c_{n-2,n-2}c_{n-1,n-1}|$ solutions.
Continue this process, $(a_1,\dots, a_{n-1})$ has $|c_{1,1} \dots c_{n-1,n-1}|$ solutions.
On the other hand $|c_{1,1} \dots c_{n-1,n-1}|$ is equal to $|d_n|$.
Similarly, the number of solutions with $a_\nu=1$ is $d_\nu$.
By Lemma \ref{freeness}, the action of $\zeta_{|d'_\nu|}$ is free.
Therefore the number of solution is $c=\gcd(d_1,\dots, d_n)$.

Next we have to check the solutions $\widetilde{\mathcal{D}(N)} \cap \mathbb{C}P(d')$ are actually ideal points.
Let $a=[t;a_1; \dots ; a_n]$ be a point of $\widetilde{\mathcal{D}(N)} \cap \mathbb{C}P(d')$.
To check $a$ is an ideal point, we have to show that $a$ is on the closure of 
$\phi(\mathcal{D}(N))$ in $\mathbb{C}P(-1,d')$. 
Because $\widetilde{\mathcal{D}(N)}$ is defined by $n-1$ equations, $\mathrm{dim} (\widetilde{\mathcal{D}(N)}) \geq 1$.
If $\widetilde{\mathcal{D}(N)}$ and $\mathbb{C}P(d')$ intersects transversely at $a$, the solutions near $a$ 
are included in the complement of the hyperplane $\mathbb{C}P(d')$ 
in $\mathbb{C}P(-1,d')$ i.e. the solutions are included in 
$\phi(\mathcal{D}(N))=\widetilde{\mathcal{D}(N)} -\mathbb{C}P(d')$ near $a$.
So we will show that $\widetilde{\mathcal{D}(N)}$ and $\mathbb{C}P(d')$ intersect at $a$ transversely.
We embed $a$ in $U_n$ as above.
Lift $a \in U_n$ to $\mathbb{C}^n$, we consider the system of equation
\begin{equation}
\label{neighborhood}
R_i
=\pm \displaystyle\prod_{\nu =1}^{n-1} a_{\nu}^{r(I)_{i,\nu}} (1-a_{\nu} t^{d_{\nu}})^{\overline{r(I)}_{i,\nu}}=  1
\quad (i=1, \dots n-1).
\end{equation}
This is the equations (\ref{replaced}) substituted $a_n$ by $1$.
Take a $\log$ of $R_i = 1$, 
\[
\log (R_i)=
\displaystyle\sum_{\nu=1}^{n-1} \biggl( r(I)_{i,\nu} \log(a_{\nu}) + \overline{r(I)_{i,\nu}}\log (1-a_\nu t^{d'_\nu}) 
\biggr) = k \pi \sqrt{-1}
\]
for some integer $k$.
Because we are interested in local behavior, 
we do not have to consider the ramification of the logarithm.
The Jacobian at $(0,a_1, \dots, a_n)$ is
\[
J(\log (R))=
\begin{pmatrix}
y_1 & \frac{r(I)_{1,1}}{a_1} & \hdots & \frac{r(I)_{1,n-1}}{a_n} \\
\vdots & \vdots & & \vdots \\
y_{n-1} & \frac{r(I)_{n-1,1}}{a_1} & \hdots & \frac{r(I)_{n-1,n-1}}{a_n} 
 \end{pmatrix}
= \begin{pmatrix} y & A \end{pmatrix}
\]
where $y_i=\displaystyle\sum_{\nu: d'_{\nu}=1} -\overline{r(I)_{i,\nu}} a_\nu $
and $A$ is the remaining $(n-1) \times (n-1)$ matrix.
Because the determinant of $A$ is a nonzero multiple of $d_n$, the determinant of $A$ is nonzero.
Therefore $A$ is invertible and the rank of $A$ is $n-1$.
Let $s \frac{\partial}{\partial t} + v_1 \frac{\partial}{\partial a_1} + 
\dots v_{n-1} \frac{\partial}{\partial a_{n-1}} \in T_a \widetilde{\mathcal{D}(N)}$.
We denote $v=(v_1, \dots , v_{n-1} )^T$.
Because $J(\log (R)) \begin{pmatrix} s \\ v \end{pmatrix} = 0$,
we have $sy+Av=0$ and $v=-sA^{-1}y$. 
$\begin{pmatrix} s \\ v \end{pmatrix} \in T_a \mathbb{C}P(d')$ means $s=0$.
Therefore, if $\begin{pmatrix} s \\ v \end{pmatrix} \in T_a \widetilde{\mathcal{D}(N)} \cap T_a \mathbb{C}P(d')$,
$\begin{pmatrix} s \\ v \end{pmatrix}=0$.
This means that $T_a \widetilde{\mathcal{D}(N)} \cap T_a \mathbb{C}P(d') = \{ 0 \}$ i.e. 
$T_a \widetilde{\mathcal{D}(N)}$ and $T_a \mathbb{C}P(d')$ intersect transversely.
\end{proof}
The above proof also shows 
\begin{corollary}
\label{smoothness}
A neighborhood of each ideal point corresponding to $(d_1, \dots, d_n)$ ($d_{\nu} > 0$ or $d_{\nu} <0$)
is smoothly embedded in $\mathbb{C}P(d'_1, \dots, d'_n,-1)$.
\end{corollary}

\begin{remark}
We can easily compute the order of $M$ and $L$ at each ideal point corresponding to $d_1, \dots ,d_n$.
In fact, we have
\[
\begin{split}
v(M) & =\displaystyle\sum_{\nu} \biggl( m'_{\nu} v(z_\nu) + m''_\nu v(w_\nu) \biggr) 
 = m \wedge (|d'_1| \rho_{i_1}, \dots, |d'_n| \rho_{i_n}), \\
v(L) & = l \wedge (|d'_1| \rho_{i_1}, \dots, |d'_n| \rho_{i_n}).
\end{split}
\]
\end{remark}

\begin{remark}
Segerman studied ideal points of once-punctured torus bundles over $S^1$ in \cite{segerman}.
He constructed ideal points of the deformation variety corresponding to any incompressible surfaces except for the fiber 
and semi-fiber case.
In the paper, it is likely to be more difficult to detect an ideal point when some $d_{\nu}=0$.
\end{remark}

\subsection{Ideal points of $X(N)$}
\label{psl-ideal}
Let $\rho \in R(N)$.
We define a regular function $I_\gamma \in \mathbb{C}[R(N)]$ by $I_\gamma(\rho)=(\mathrm{tr}(\rho))^2$.
Since $I_\gamma$ is invariant under the action of $\PSLC$,
this is also a regular function of $X(N)$.
Let $z=(z_1, \dots, z_n) \in \mathcal{D}(N)$.
We denote the corresponding representation by $\rho_z$ (this is well-defined up to conjugation).
By conjugation, we assume that $\rho_z(\mathcal{M})$ and $\rho_z(\mathcal{L})$ are 
diagonal matrices:
\[
\rho_z(\mathcal{M})=\begin{pmatrix} \mu & 0 \\ 0 & \mu^{-1} \end{pmatrix}, \quad 
\rho_z(\mathcal{L})=\begin{pmatrix} \lambda & 0 \\ 0 & \lambda^{-1} \end{pmatrix}. 
\]
Then we have $M=\mu^2$ and $L=\lambda^2$ (see \cite{neumann-zagier}).
$I_{\mathcal{M}}=(\mu+\mu^{-1})^2=\mu^2+\mu^{-2}+2=M+M^{-1}+2$ is a regular function of $\mathcal{D}(N)$ and also $X(N)$.
We have $v(I_{\mathcal{M}})=\min(v(M),-v(M))=-|v(M)|$.
For $\rho_z$, we have 
\[
v(I_\mathcal{M})=-|v(M)|=-| m \wedge (|d'_1| \rho_{i_1}, \dots, |d'_n| \rho_{i_n})|=-|m \wedge x|, \quad 
v(I_\mathcal{L})=-|l \wedge x|
\]
where we denote $x=(|d'_1| \rho_{i_1}, \dots, |d'_n| \rho_{i_n})\in\mathbb{Z}^{2n}$.
If $m \wedge x$ or $l \wedge x$ is nonzero, $v$ corresponds to an ideal point of 
$X(N)$ because the character of $\mathcal{M}$ or $\mathcal{L}$ diverges.
If $v(\mathcal{M}^p \mathcal{L}^q)=0$, $\mathcal{M}^p \mathcal{L}^q$ represents the boundary slope of the ideal
point (see \cite{ccgls}).
Because $ ( (- l \wedge x) m + (m \wedge x) l ) \wedge x = 0$, 
$\mathcal{M}^{- l \wedge x} \mathcal{L}^{m \wedge x}$ represents the boundary slope.
We summarize our results:
\begin{theorem}
Let $N$ be an orientable 3-manifold with torus boundary.
Let $\mathcal{M}$ and $\mathcal{L}$ be generators of $\pi_1(\partial N)$.
Let $K$ be an ideal triangulation of $N$ such that $\mathcal{D}(N,K)$ is non-empty.
For a degenerate index $I$ which satisfies that $d(I)>0$ or $d(I)<0$ and  
$m \wedge (d_1 \rho_{i_1}, \dots, d_n \rho_{i_n})$ or $l \wedge (d_1 \rho_{i_1}, \dots, d_n \rho_{i_n})$
is nonzero, there exists an incompressible surface with the boundary slope 
$\mathcal{M}^{(-l \wedge (d_1 \rho_{i_1}, \dots, d_n \rho_{i_n}))}
\mathcal{L}^{( m \wedge (d_1 \rho_{i_1}, \dots, d_n \rho_{i_n}))}$.
\end{theorem}

\section{Examples}
\label{sec:examples}
We give some examples from SnapPea's census manifolds \cite{weeks}.
I used Snap \cite{snap} to show gluing equations.
We remark that a census triangulation has non-empty solution $\mathcal{D}(N)$ 
because there is a discrete faithful character.
\subsection{m006}
m006 has 3 ideal simplices.
The gluing equations are given at 3 edges by
\[
\begin{split}
z_1^2 w_1^{-1} z_2 w_2 z_3 w_3 =1 \\
w_1^{-1} z_2^{-1} z_3^{-1}=1 \\
z_1^{-2} w_1^2 w_2^{-1} w_3^{-1} =1. 
\end{split}
\]
One of the equations is dependent on the others,
so we only consider the last two equations.
The matrix $R$ corresponding to (\ref{R}) is 
\[
R = 
\begin{pmatrix} 
0 & -1 & -1 & 0 & -1 & 0 \\
-2 & 2 & 0 & -1 & 0 & -1 
\end{pmatrix}.
\]
$m$ and $l$ are given by
\[
m=(0,0,1,0,-1,0), \quad l=(0,1,0,1,2,-1).
\]
Consider the $3^3$ combinations of degeneration indices, the indices satisfying 
$d(I) >  0$ or $d(I) < 0$ are 
\begin{equation}
\label{m006}
(1,1,\infty ), \quad
(1,\infty , 1), \quad
(0,0,\infty ) , \quad
(0,\infty, 0).
\end{equation}
The degeneration matrix corresponding to the first degeneration index is 
\[
R( ( 1,1,\infty ) )=
\begin{pmatrix}
-1 & 0 & 1\\
2 & -1 & 1
\end{pmatrix}.
\]
So the degeneration vector is given by 
\[
\left(\det \begin{pmatrix} 0 & 1 \\ -1 & 1 \end{pmatrix},-\det \begin{pmatrix} -1 & 1 \\ 2 & 1\end{pmatrix},
\det \begin{pmatrix} -1 & 0 \\ 2 & -1 \end{pmatrix} \right) = (1,3,1).
\]
The order of $M$ at the ideal point is 
\[
v(M)=m \wedge (1 \cdot \rho_1, 3 \cdot \rho_1, 1 \cdot \rho_{\infty})=m \wedge (1,0,3,0,-1,1)=-1.
\]
Similarly, the degeneration vectors corresponding to (\ref{m006}) are
\[
(1,3,1), \quad (-1,-1,-3), \quad (-1,-2,-2), \quad (1,2,2) 
\]
respectively.
These four degeneration vectors correspond to ideal points by Theorem \ref{main-theorem}.
The orders $(v(M),v(L))$ of $M$ and $L$ at each ideal point are 
\[
(-1,-3), \quad (1,3), \quad (-4,2), \quad (4,-2)
\]
respectively.
The corresponding boundary slopes are $-3/1$, $-3/1$, $1/2$, $1/2$ respectively.

The degeneration index $(\infty ,0,0)$ is a typical example of subsection \ref{compute-cone}.
For this index, we have the degeneration vector $(0,0,0)$, the degeneration matrix is 
\[
R((\infty,0,0))=
\begin{pmatrix} 1 & -1 & -1 \\ 0 & 0 & 0 \end{pmatrix}.
\]
For $(\epsilon_1, \epsilon_2, \epsilon_3) =(0,1,1)$, $(1,0,1)$ and $(1,1,0)$, 
we have the degeneration vector $(-1,1)$, $(1,1)$ and $(1,1)$ respectively.
Therefore $(\rho_\infty,0,\rho_0)$ and $(\rho_\infty, \rho_0, 0)$ are generators of $H(\infty,0,0) \cap [R]^{\bot}$.

\subsection{m009}
This manifold is homeomorphic to the once-punctured torus bundle with monodromy 
\[
\phi=  \begin{pmatrix} 1 & 0 \\ 2 & 1 \end{pmatrix} \begin{pmatrix} 1 & 1 \\ 0 & 1 \end{pmatrix}
= \begin{pmatrix} 1 & 1 \\ 2 & 3 \end{pmatrix}.
\]
The gluing equations are given by 
\[
z_1^2 z_2^2 z_3^2 =1, \quad w_1^{-1} z_2^{-2} w_2^2 w_3^{-1} =1 .
\]
The matrix $R$ is 
\[
R = 
\begin{pmatrix} 
2 & 0 & 2 & 0 & 2 & 0 \\
0 & -1 & -2 & 2 & 0 & -1 
\end{pmatrix}.
\]
$m$ and $l$ are given by
\[
m = (-1, -1, -2, 1, -1, 0), \quad l = (3, -1, 2, -1, -1, 2).
\]
Consider the $3^3$ combinations of degeneration indices, the indices satisfying 
$d(I) >  0$ or $d(I) < 0$ are 
\begin{equation}
\label{m009}
(0,0,\infty ), \quad 
(\infty ,0,0).
\end{equation}
The degeneration vectors corresponding to (\ref{m009}) are
\[
(-2,-2,-4), \quad (4,2,2),
\]
respectively.
By Theorem \ref{main-theorem}, each of the two degeneration vectors has corresponding two ideal points.
The orders $(v(M),v(L))$ of $M$ and $L$ at the ideal points corresponding to $(-2,-2,-4)$ are $(1,-3)$ and
those ideal points corresponding to $(4,2,2)$ are $(-1,3)$.
The boundary slopes corresponding to these ideal points are $3/1$.
We will study these two ideal points.

By eliminating $z_1$, the edge relations can be reduced to the following equation of 2 variables:
\begin{equation}
\label{factorization}
\begin{split}
& z_2^2 (1-z_3)^2 - z_3^2 (z_2^2(1-z_3)-(1-z_2)^2)^2 \\
& =(z_2(1-z_3)+z_3(z_2^2(1-z_3)-(1-z_2)^2)) \cdot \\
& \quad \quad \quad \quad (z_2(1-z_3)-z_3(z_2^2(1-z_3)-(1-z_2)^2))=0.
\end{split}
\end{equation}
By using Snap, we can see that the complex parameters corresponding to a discrete faithful representation are given by 
$z_1=z_3=\displaystyle\frac{1+\sqrt{7}i}{2}, z_2=\displaystyle\frac{3+\sqrt{7}i}{8}$.
By substituting these variables, we can show that
the first factor of (\ref{factorization}) contains the discrete faithful representation.
By comparison with the $\mathrm{SL}(2,\mathbb{C})$ A-polynomial, the second factor does not lift to 
any component of the $\mathrm{SL}(2,\mathbb{C})$ representations. 

Consider the ideal points corresponding to $I=(\infty,0,0)$.
Since $d(I)=(2,2,4)$, we have $d'(I)=(1,1,2)$.
Let $z_1=\frac{1}{a_1 t^2}$, $z_2=a_2 t$, $z_3=a_3 t$.
Then the equations (\ref{m009}) are replaced by 
\begin{equation}
\frac{a_2^2 a_3^2}{a_1^2}=1, \quad -\frac{a_1 (1-a_2 t)^2}{(1-a_1 t^2) a_2^2 (1-a_3 t)}=1.
\end{equation}
Considering these equations in $U_3$ (i.e. $a_3=1$) and at infinity (i.e. $t=0$), we have
\[
a_1^{-2}a_2^2=1, \quad -a_1 a_2^{-2}=1. 
\] 
The solutions are $(a_1,a_2)=(-1,1),(-1,-1)$.
For the first solution $(a_1,a_2)=(-1,1)$, 
we have $z_1=\frac{1}{-t^2}+o(t^{-1})$, $z_2= t+o(t^2)$, $z_3=t+o(t^2)$.
Since the ideal point is smooth by Corollary \ref{smoothness}, these are convergent power series.
Substitute these variables, then the factors of (\ref{factorization}) are
\[
\begin{split}
z_2(1-z_3) \pm z_3(z_2^2(1-z_3)-(1-z_2)^2) & = t \cdot 1 + t \cdot (t^2 - 1)+(\textrm{higher order}) \\
& = t \pm (-t)+(\textrm{higher order}).
\end{split}
\]
So we can conclude that the ideal point corresponding to $(a_1,a_2,a_3)=(-1,1,1)$ is an ideal point of 
the first factor of (\ref{factorization}).
From the second solution $(a_1,a_2,a_3)=(-1,-1,1)$, we can conclude that the corresponding ideal point is contained 
in the second factor of (\ref{factorization}).

Finally we remark that m009 has another ideal point which can not be detected by our method.
Consider the degeneration index $I_1=(0,1,1)$.
Then the degeneration vector $d(I_1)$ is equal to $(0,2,4)$.
This satisfies the necessary condition $d(I_1) \geq 0$ but does not satisfy our sufficient condition.
But this suggests that we should put $z_2=1-a_2 t$ and $z_3=1-a_3 t^2$ and consider gluing equations 
in $\mathbb{C}P(-1,1,2) \times (\mathbb{C}-\{0,1\})$ with the coordinate $([t;a_2;a_3],z_1)$.
The gluing equations at infinity are
\[
z_1^2=1, \quad (1-z_1)^{-1} a_2^2 a_3^{-1}=1,
\]
so we  have a solution at infinity $([0; 1 ; 2],-1)$.
Next we consider the degeneration index $I_2=(0,\infty,1)$.
We have $d(I_2)=(2,2,0)$.
So $d(I_2)$ is also a degeneration vector which satisfies the necessary condition but does not satisfy our condition. 
As in the case of $I_1$ we should put $z_1=a_1 t$ and $1/z_2= a_2 t$ and consider gluing equations 
in $\mathbb{C}P(-1,1,1) \times (\mathbb{C}-\{0,1\})$ with the coordinate $([t;a_1;a_2],z_3)$. 
But in this case, the gluing equations at infinity are 
\[
a_1^2 a_2^{-2} z_3^2 = 1, \quad (1-z_3)^{-1}=1.
\]
This has no solution in $\mathbb{C}P(-1,1,1) \times (\mathbb{C}- \{0,1\})$ because $z_3=0$ by the 
second equation.
These two examples show that it is subtle problem to find ideal points of $\mathcal{D}(N)$ 
when some ideal tetrahedron does not degenerate.



\begin{thebibliography}{99}
\bibitem{boyer-zhang}
S. Boyer, X. Zhang, {On Culler-Shalen seminorms and Dehn filling}, 
Ann. of Math. (2) 148 (1998), no. 3, 737--801.


\bibitem{ccgls}
D. Cooper, M. Culler, H. Gillet, D. Long, P. Shalen,
{\it Plane curves associated to character varieties of $3$-manifolds},
Invent. Math. 118 (1994), no. 1, 47--84. 


\bibitem{culler-shalen}
M. Culler, P. Shalen, {\it Varieties of group representations and splittings 
of $3$-manifolds}, Ann. of Math. (2) 117 (1983), no. 1, 109--146. 



\bibitem{neumann}
W.D. Neumann, {\it Combinatorics of triangulations and the
Chern--Simons invariant for hyperbolic $3$--manifolds}, from: ``Topology 90,
Proceedings of the Research Semester in Low Dimensional Topology at
Ohio State'', Walter de Gruyter Verlag, Berlin--New York (1992)
243--272


\bibitem{snap}
O. Goodman, {\it Snap} , computer program.


\bibitem{heusener-porti}
M. Heusener, J. Porti, 
{\it The variety of characters in ${\rm PSL}\sb 2(\mathbb C)$},
Bol. Soc. Mat. Mexicana (3) 10 (2004), Special Issue, 221--237. 

\bibitem{neumann-zagier}
W.D.Neumann, D.Zagier,  {\it Volumes of hyperbolic 3-manifolds},
Topology 24 (1985), 307-332.

\bibitem{segerman}
H. Segerman,
{\it Incompressible Surfaces in Hyperbolic Punctured Torus Bundles are Strongly Detected},
math.GT/0610302.


\bibitem{tillmann-ideal}
S. Tillmann, {\it Boundary slopes and the logarithmic limit set},
Topology 44 (2005), no. 1, 203--216.

\bibitem{tillmann-normal}
S. Tillmann, {\it Normal surfaces in topologically finite 3-manifolds},
math.GT/0406271.

\bibitem{yoshida}
T. Yoshida, {\it On ideal points of deformation curves of hyperbolic $3$-manifolds with one cusp},
Topology 30 (1991), no. 2, 155--170.

\bibitem{weeks}
J. Weeks, {\it SnapPea}, computer program. 

\end{thebibliography}
\end{document}